\documentclass[12pt]{article}
\usepackage{amsmath}
\usepackage{amssymb}
\usepackage{latexsym}
\usepackage{amsfonts}
\date{}

\newcommand{\qed}{\ifmmode$\Box$\else{\unskip\nobreak\hfil
\penalty50\hskip1em\null\nobreak\hfil$\Box$
\parfillskip=0pt\finalhyphendemerits=0\endgraf}\fi}

\newtheorem{Pa}{Paper}[section]
\newtheorem{Tm}[Pa]{{\bf Theorem}}

\newtheorem{Cy}[Pa]{{\bf Corollary}}

\newtheorem{Pn}[Pa]{{\bf Proposition}}

\newtheorem{Ee}[Pa]{{\bf Example}}
\newtheorem{Dn}[Pa]{{\bf Definition}}

\newcommand{\CC}{{\mathchoice
{\setbox0=\hbox{$\displaystyle\rm C$}\hbox{\hbox
to0pt{\kern0.4\wd0\vrule height0.9\ht0\hss}\box0}}
{\setbox0=\hbox{$\textstyle\rm C$}\hbox{\hbox
to0pt{\kern0.4\wd0\vrule height0.9\ht0\hss}\box0}}
{\setbox0=\hbox{$\scriptstyle\rm C$}\hbox{\hbox
to0pt{\kern0.4\wd0\vrule height0.9\ht0\hss}\box0}}
{\setbox0=\hbox{$\scriptscriptstyle\rm C$}\hbox{\hbox
to0pt{\kern0.4\wd0\vrule height0.9\ht0\hss}\box0}}}}

\begin{document}
\title{Rational functions associated with the white noise space
and related topics}
\author{Daniel Alpay and David Levanony}
\maketitle

\begin{abstract}
Motivated by the hyper-holomorphic case we introduce and study
rational functions in the setting of Hida's white noise space.
The Fueter polynomials are replaced by a basis computed in terms
of the Hermite functions, and the Cauchy-Kovalevskaya product is
replaced by the Wick product.
\end{abstract}




\section{Introduction}
\setcounter{equation}{0}
Consider a $K$-vector space $V$ spanned by a family of functions
$(f_\alpha)_{\alpha\in\ell}$, where $\ell$ is a countable set of
indices, and where $K$ denotes either ${\mathbb R}$, ${\mathbb
C}$, or the skew field of quaternions ${\mathbb H}$ (in this latter
case we assume that $V$ is a right vector space, to fix the
ideas). Define on $V$ an operation by
\begin{equation}
\label{law} f_\alpha\circ f_\beta=f_{\alpha+\beta},\quad
\alpha,\beta\in\ell.
\end{equation}
In general, such a product depends on the basis and need not carry
any structure related to $V$. There are at least two cases we are
aware of, where the multiplication-like law \eqref{law} carries
much information. The first is the case of hyper-holomorphic
functions. This case corresponds to $K={\mathbb H}$ and
$\ell={\mathbb N}^3$. The $(f_\alpha)$ are the Fueter monomials
and $\circ$ is the Cauchy-Kovalevskaya product; see \cite{bds},
\cite{MR618518}. The second case, which is the topic of the
present work, corresponds to $K={\mathbb R}$, and $\ell$ the
space of sequences
$(\alpha_n)_{n\in{\mathbb N}}$ of integers for which $\alpha_n=0$
for $n$ sufficiently large; the space $V$ is Kondratiev's space
(which includes Hida's white noise space and Hida's space of
distributions), with the $(f_\alpha)$ forming an orthonormal
Hilbert space basis of the white noise space built in terms of the
Hermite functions, and $\circ$ is the Wick product. All these
notions will be
reviewed in the sequel.\\

The notions of rational functions, de Branges Rovnyak spaces and
Schur-Agler classes were introduced for the case of
hyper-holomorphic functions in the papers \cite{assv}, \cite{as1},
\cite{asv-cras}, \cite{MR2124899}, \cite{MR2240272}. The case of
Clifford algebra valued functions has been considered in
\cite{MR2275397}. In the approach developed in these papers,
important tools were the study of the Gleason problem and the
introduction of counterparts of the Leibenzon operators for
hyper-holomorphic functions. The purpose of this paper is to make
a similar study within the white noise space setting. Results obtained
in this paper are to be applied to problems in stochastic system
theory, a work to be summarized in a future publication.\\

To provide further motivation, it is best to first take a detour via
several complex variables, and discuss Gleason's problem and the
Leibenzon operators. Recall that the backward shift operator
\[R_0f(z)=\frac{f(z)-f(0)}{z}
\]
plays an important role in operator theory and in the theory of
linear systems; see \cite{Rudin}, \cite{bgk1}. It has
counterparts in several complex variables, as we now recall: Let
$f$ be a function of $N$ complex variables, analytic in a
neighborhood of the origin. Then (see e.g. \cite[p.
151]{MR54:11066}), it holds that
\begin{equation}
\label{bw2}
f(z)-f(0)=\int_0^1\frac{\rm d}{\rm dt}f(tz)dt
         =\sum_{j=1}^Nz_j ({\mathcal R}_jf)(z)
\end{equation}
where ${\mathcal R}_j$ denotes the Leibenzon's backward shift
operator (see \cite[p. 117-118]{rudin-ball})
\begin{equation}
\label{bastille}
{\mathcal R}_jf(z)=\int_0^1\frac{\partial
f(tz)}{\partial z_j}dt= \sum_{\alpha\in{\mathbb
N}^N}\frac{\alpha_j}{|\alpha|}c_\alpha z^{\alpha-{\mathbf u}_j},
\end{equation}
with
\begin{equation}
\label{bastille2} f(z)=\sum_{\alpha\in{\mathbb N}^N}c_\alpha
z^\alpha.
\end{equation}
In \eqref{bastille}--\eqref{bastille2}, we have used the the
multi-index notation, and, for $j=1,\ldots , N$, the symbol
${\mathbf u}_j$ denotes the
index with all components equal to $0$, with the exception of the
$k$-th, equal to $1$.
Moreover, $\alpha-{\mathbf u}_j$ is defined to be $(0,0,\ldots, 0)$ when
one of its entries is strictly negative.\\

Comparing the definition of $R_0$ with \eqref{bw2}, suggests that
the operators ${\mathcal R}_j$ are a generalization of the
backward--shift operator $R_0$. This is indeed the case, but the
situation is more complex: Gleason's problem for a space
${\mathcal M}$ of functions analytic in a neighborhood of the
origin, asks the following: Can we write for $f\in{\mathcal M}$
\[
f(z)-f(0)=\sum_{j=1}^Nz_jf_j(z)\]
where $f_j\in{\mathcal M}$? The operators ${\mathcal R}_j$ allow
to solve Gleason's problem in various spaces of power series. We
refer to \cite{ad4} for a study of these operators and to
\cite{adr1} for an extension of the Beurling--Lax theorem using
Gleason's problem in the setting of the ball.\\

The paper is written with three different audiences in mind,
namely researchers in multi-dimensional system theory,
hyper-holomorphic functions and stochastic analysis. It consists
of seven sections besides the introduction and is organized as
follows: In Section 2 we briefly review the hyper-holomorphic
case. This is important because we want to specify the
similarities between the Cauchy-Kovalevskaya product and the Wick
product. In Section 3 we review the main properties of Hida's
white noise space. Rational functions in the white noise space
are defined and studied in Section 4. In Sections 5 and 6 we study
several counterparts of classical spaces in the stochastic
setting. Section 5 is devoted to the Arveson space and its
multipliers, while Section 6 is devoted to the counterpart of the
Schur-Agler classes of the polydisk. In the last section we prove a
uniqueness result related to the Leibenzon operators in certain
spaces of power series in a
countable number of variables.\\

A table showing the parallels between the hyper-holomorphic case
and the case of the white noise space is provided at the end of
the paper.\\


\section{The hyper-holomorphic case: a short review}
We review the main features of the hyper-holomorphic case relevant
to the present paper. First recall that the skew-filed of
quaternions consists of elements of the form
\[x=x_0{\mathbf e_0}+x_1{\mathbf e_1}+x_2{\mathbf e_2}+x_3{\mathbf
  e_3},
\]
where the $x_i\in{\mathbb R}$ and the ${\mathbf e_i}$ satisfy the
Cayley multiplication table
\begin{equation*}\begin{array}{|c|c|c|c|c|} \hline
   &\bf e_0&\bf e_1&\bf e_2&\bf e_3\\ \hline
 \bf e_0&\bf e_0 &\bf e_1&\bf e_2&\bf e_3\\ \hline \bf e_1&\bf
e_1&-\bf e_0&\bf e_3&-\bf e_2\\ \hline \bf e_2&\bf e_2&-\bf e_3&-\bf
e_0&\bf e_1\\ \hline \bf e_3&\bf e_3&\bf e_2&-\bf e_1&-\bf e_0 \\
\hline \end{array} \end{equation*}
We set ${\mathbf e_0}=1$.\\

The function $f:\Omega\subset{\mathbb
R}^4\rightarrow{\mathbb H}$ is called {\sl left hyper-holomorphic}
if
\begin{equation}
\label{Stalingrad-sur-la-ligne-de-metro-numero-7}
D f{:=} \frac{\partial}{\partial x_0}f+{\mathbf e}_{\mathbf 1}
\frac{\partial}{\partial x_1}f+ {\mathbf e}_{\mathbf 2}
\frac{\partial}{\partial x_2}f+{\mathbf e}_{\mathbf 3}
\frac{\partial}{\partial x_3}f=0.
\end{equation}
Write $f=f_0+{\mathbf e_1}f_1+{\mathbf e_2}f_2+{\mathbf e_3}f_3$.
The components $f_j$ of $f$ satisfy the system of equations

\begin{equation}
\begin{split}
\frac{\partial f_0}{\partial x_0} -\frac{\partial f_1}{\partial
x_1}-\frac{\partial f_2}{\partial x_2}-
\frac{\partial f_3}{\partial x_3}&=0,\\
\frac{\partial f_0}{\partial x_1} +\frac{\partial f_1}{\partial
x_0}-\frac{\partial f_2}{\partial x_3}+
\frac{\partial f_3}{\partial x_2}&=0,\\
\frac{\partial f_0}{\partial x_2} +\frac{\partial f_1}{\partial
x_3}+\frac{\partial f_2}{\partial x_0}-
\frac{\partial f_3}{\partial x_1}&=0,\\
\frac{\partial f_0}{\partial x_3} -\frac{\partial f_1}{\partial
x_2}+\frac{\partial f_2}{\partial x_1}+ \frac{\partial
f_3}{\partial x_0}&=0.
\end{split}
\label{23-juillet-2002}
\end{equation}
One can apply to this system of partial differential equations
the Cauchy--Kovalevskaya theorem: Let $\varphi(x_1,x_2,x_3)$ be a
real analytic function from some open domain of ${\mathbb R}^3$
into ${\mathbb H}$, that is, $\varphi$ is given by four coordinate
real-analytic, real--valued functions
$$
\varphi(x_1,x_2,x_3)=\varphi_0(x_1,x_2,x_3)+\sum_1^3{\mathbf e_i}
\varphi_i(x_1,x_2,x_3).$$
The Cauchy--Kovalevskaya theorem, see \cite[Section
1.7]{MR93j:26013}, implies that the system of equations
\eqref{23-juillet-2002}, with initial conditions
$$f_i(0,x_1,x_2,x_3)=\varphi_i(x_1,x_2,x_3)$$
admits a unique real analytic solution in a neighborhood of the
origin in ${\mathbb R}^4$. This solution
$$f(x_0,x_1,x_2,x_3)=f_0(x_0,x_1,x_2,x_3)+
\sum_1^3{\mathbf e_i}f_i(x_0,x_1,x_2,x_3)$$ is hyper-holomorphic
by definition and is called the Cauchy--Kovalevskaya extension of
the function $\varphi$. We will use the notation
\[
f=\mathbf{CK}(\varphi).
\]
Let $\alpha=(\alpha_1,\alpha_2,\alpha_3)\in{\mathbb N}^3$ and let
$\varphi(x)=x_1^{\alpha_1}x_2^{\alpha_2}x_3^{\alpha_3}\stackrel{\rm
def.}{=} x^\alpha$. The corresponding hyper-holomorphic function
is the Fueter monomial $\zeta^\alpha$. The case where
$\varphi(x)=x_\ell$ ($\ell=1,2,3$) leads to the hyper-holomorphic
variables
\[
\zeta_\ell(x)=x_\ell-{\mathbf e}_\ell x_0,\quad \ell=1,2,3,
\]
that is $\zeta_\ell=\mathbf{CK}(x_\ell)$. The notations
\begin{equation}
\label{zeta}
\zeta(x)=\begin{pmatrix}\zeta_1(x)&\zeta_2(x)&\zeta_3(x)\end{pmatrix}
\end{equation}
and
\begin{equation}
\label{zetaN} \zeta^{(N)}(x)=\begin{pmatrix}\zeta_1(x) I_N
&\zeta_2(x) I_N &\zeta_3(x)I_N\end{pmatrix}
\end{equation}
will prove useful.\\


The point-wise product of two hyper-holomorphic functions, say
$f$ and $g$, need not be hyper-holomorphic. Their
Cauchy-Kovalevskaya product $f\circ g$ has been introduced in
1981 by F. Sommen in \cite{MR618518} and is defined as the
Cauchy-Kovalevskaya extension of
\[
f(0,x_1,x_2,x_3)g(0,x_1,x_2,x_3).
\]
It is a hyper-holomorphic function. We set ${\bf R}$ to be the
restriction of a hyper-holomorphic function to $x_0=0$. Then,
with $x=\begin{pmatrix} x_1&x_2&x_3\end{pmatrix}$ we have
\[
f\circ g=\mathbf{CK}({\mathbf R}(f){\mathbf R}(g)).
\]
Hence, we have that $f\circ g=g\circ  f$ if and only if the
quaternionic-valued functions $\mathbf R(f)$ and
$\mathbf R(g)$ commute.\\

Note that for the Fueter monomials we have
\[
\zeta^\alpha\circ\zeta^\beta=\zeta^{\alpha+\beta},\quad
\alpha,\beta\in{\mathbb N}^3.
\]

Every function $f$ hyper-holomorphic in a neighborhood of the
origin can be written as a power series expansion using the Fueter
monomials $\zeta^\alpha$
\[
f=\sum_{\alpha\in{\mathbb N}^3}\zeta^\alpha f_\alpha,\quad
f_\alpha\in{\mathbb H},
\]
and the Cauchy-Kovalevskaya product has a nice interpretation in
terms of these expansions: It is a convolution product, also
called the Cauchy product:
\[
f\circ g=\sum_{\alpha\in{\mathbb
N}^3}\zeta^\alpha\left(\sum_{\beta\le \alpha} f_\beta
g_{\alpha-\beta}\right).
\]
To define rational functions we first need another definition and a
result: If $g$ is a ${\mathbb H}^{N\times N}$--valued
hyper-holomorphic function, we denote
\[
g^{n\circ}=g\circ\cdots \circ g \quad (n\quad{\rm times}).
\]
If moreover, $g(0)=0$, then the series
\begin{equation}
\label{inv}
(I_N-g)^{-\circ}\stackrel{\rm
def.}{=}\sum_{n=0}^\infty g^{n\circ}
\end{equation}
converges in a neighborhood of the origin to a hyper-holomorphic function.\\

\begin{Tm}
\label{110108} The following conditions are equivalent:
\begin{enumerate}
\item
The function $f(0,x_1,x_2,x_3)$ is a rational function of the
three real variables $x_1,x_2$ and $x_3$ with values in ${\mathbb
H}$ (that is, each of its real components is a rational function
of the three real variables $x_1,x_2$ and $x_3$ ),
and analytic
at the origin.

\item We can write
\begin{equation}
\label{real}
f(x)=D+C\circ(I_N-\zeta^{(N)}(x) A)^{-\circ}\circ
\zeta^{(N)}(x) B
\end{equation}
where $A,B,C$ and $D$ are matrices with entries in ${\mathbb H}$
of appropriate dimensions.

\item $f$ is obtained from the Fueter monomials and the
quaternions after a finite number of additions,
Cauchy-Kovalevskaya multiplications and inversions (the latter
defined by \eqref{inv}).
\end{enumerate}
\label{tototototot}
\end{Tm}

Expression \eqref{real} is called a {\sl realization of $f$
centered at the origin}, and comes from system theory; see
\cite{bgk1}, \cite{MR54:2342}, \cite{MR80c:93028}. Using the
abuse of notation
\[
\zeta^{(N)}(x)=\zeta,
\]
we will write this expression as
\[
f(x)=D+C\circ(I-\zeta A)^{-\circ}\circ \zeta B.
\]
More explicitly, one can also write \eqref{real} as
\[
f(x)=D+C(I_N-\zeta_1A_1-\zeta_2A_2-\zeta_3A_3)^{-\circ}\circ
(\zeta_1B_1+\zeta_2B_2+\zeta_3B_3)
\]
for matrices $A_1,B_1,\cdots$ of appropriate dimensions.
Note also that
\begin{equation}
\label{real4} {\mathbf R}(D+C\circ(I-\zeta A)^{-\circ}\circ \zeta
B)=D+C(I-xA)^{-1}xB,
\end{equation}
and
\begin{equation}
\label{real5} \mathbf{CK}(D+C(I-xA)^{-1}xB)=
D+C\circ(I-\zeta A)^{-\circ}\circ \zeta B.
\end{equation}
\begin{Dn}
A ${\mathbb H}^{p\times q}$-valued function $f$ hyper-holomorphic
in a neighborhood of the origin, is called {\rm rational} if any
of the equivalent conditions in the previous theorem is in force.
\end{Dn}

Theorem \ref{tototototot} is proved in the above mentioned
papers. Here we will give a short and slightly different proof,
both for completeness and
because similar arguments will be used in Section \ref{Anaelle}.\\

{\bf Proof of Theorem \ref{tototototot}:} We begin by proving the
equivalence between the first two conditions. We recall that any
rational ${\mathbb C}^{p\times q}$-valued function $W(z)$ of $N$
complex variables $z_1,\ldots, z_M$ can be written as
\begin{equation}
\label{real2} W(z)=D+C(I_N-zA)^{-1}zB.
\end{equation}
In this expression, $D=W(0)\in{\mathbb C}^{p\times q}$,
$N\in{\mathbb N}$ and $C\in{\mathbb C}^{p\times N}$. Furthermore,
$zA=z_1A_1+\cdots +z_MA_M$ and $zB=z_1B_1+\cdots z_MB_M$, where
the $A_j\in{\mathbb C}^{N\times N}$ and the $B_j\in{\mathbb
C}^{N\times q}$. This result originated with \cite[Theorem 5, p.
107]{MR924203}. A different proof, based on Gleason's problem,
has been recently given in \cite{a-dubi3}. The realization
\eqref{real2} still holds for real-valued functions of three real
variables. Identifying the quaternions with ${\mathbb R}^4$, it
is also seen that the realization result still holds for
quaternionic-valued rational functions of three real variables.
Equation \eqref{real} then follows directly from \eqref{real5}.
The converse follows from
\eqref{real4}.\\

To study the equivalence with the third condition we note the
following: Since the first condition is invariant under summation,
pointwise multiplication and pointwise inversion, the sum and
product of two hyper-holomorphic functions of compatible dimensions and
the inverse (in the sense of \eqref{inv}), are still rational. The
equivalence with the third condition follows then from the fact
that a rational function of three real variables with
quaternionic entries is obtained after a finite number of sums,
products and divisions of monomials $x^\alpha$.
\mbox{}\qed\mbox{}\\

\section{The white noise space and Kondratiev's space}
In the previous section the underlying space was the space of
functions hyper-holomorphic in a neighborhood of the origin. In
the setting we will now review, the situation is more complex (at
least in the present stage of the theory). The first step is to
build the white noise space, and then to go beyond, to a space
of distributions.\\

To define Hida's white noise space first set ${\mathcal S}$ to be
the Schwartz space of rapidly decreasing functions, and ${\mathcal
S}^\prime$ its topological dual (the space of tempered
distributions). We denote by ${\mathcal F}$ the $\sigma$-algebra
of the Borel sets of ${\mathcal S}^\prime$. Hida's white noise
space is constructed as follows, using the Bochner-Minlos theorem.
First note that the function
\[K(s_1-s_2)=\exp(-\|s_1-s_2\|_{{\mathbf L}_2
({\mathbb R})}^2/2),\]
is positive in the sense of reproducing kernels
in ${\mathcal S}$. Since ${\mathcal S}$ is
nuclear, there exits a probability measure $P$ on
$({\mathcal S}^\prime, {\mathcal F})$ such that,
for all $s\in {\mathcal S}$,
\begin{equation}
\label{1967} E(e^{iQ_s})=
e^{-\frac{\|s\|_{{\mathbf L}_2({\mathbb
R})}^2}{2}},
\end{equation}
where $Q_s$ denotes the linear functional
$Q_s(s^\prime)=\langle
s^\prime,s\rangle_{{\mathcal S}^\prime,{\mathcal S}}$. See for instance
\cite[Th\'eor\`eme 2, p. 342]{MR35:7123}. Equation \eqref{1967}
implies in particular that
\begin{equation}
\label{Campo Formio} E(Q_s)=0\quad{\rm and}\quad
E(Q_s^2)=\|s\|_{{\mathbf L}_2({\mathbb R})}^2.
\end{equation}
${\mathcal W}\stackrel{\rm def.}{=} {\mathbf L}_2({\mathcal
S}^\prime, {\mathcal F},P)$ is the white noise probability space.
In accordance with the notation standard in probability theory, we
set $\Omega={\mathcal S}^\prime$. Thus,
\[
{\mathcal W}= {\mathbf L}_2(\Omega, {\mathcal F},P).\]

The white noise space ${\mathcal W}$ admits a special orthonormal
basis $(H_\alpha)$, indexed by the set $\ell$ of finite sequences
of ${\mathbb N}^{\mathbb N}$, and is built in terms of the
Hermite functions (which themselves, are constructed by the
Hermite polynomials). We refer the reader to \cite[Chapter
2]{MR1408433} and to the papers \cite{bosw}, \cite[p. 305]{eh},
where the main features of the theory are reviewed. Because of
the forthcoming definition of the Wick product (see
\eqref{la_tour_eiffel} below), it suffices to briefly recall the
definition of the $H_\alpha$ when $\alpha=e_{k}$, $k=1,2,\ldots$.
Here, we have denoted by $e_k$ the element of $\ell$ with all
entries equal to $0$, with the exception of the $k$-th, being
equal to $1$. To define $H_{e_k}$, let $\zeta_k$ be the $k$-th
Hermite function (which itself is computed in terms of the $k$-th
Hermite polynomial $h_k$). Then, $\zeta_k\in{\mathcal S}$ and
\[
H_{e_k}(\omega)=h_k(\langle \omega,\zeta_k\rangle_{\Omega,{\mathcal S}}).
\]
Every $ F\in {\mathbf L}_2(\Omega, {\mathcal F},P)$ admits a
representation
\begin{equation}
F=\sum_{\alpha\in\ell} c_\alpha H_\alpha,
\label{fps}
\end{equation}
with
\[
\sum_{\alpha\in\ell} c_\alpha^2\alpha!<\infty,
\]
called Wiener-It\^o chaos expansion; see \cite[Theorem 2.2.4 p.
23]{MR1408433}.\\

The Wick product is defined by
\begin{equation}
\label{la_tour_eiffel}
H_\alpha\lozenge H_\beta=H_{\alpha+\beta},
\end{equation}
which is reminiscent of the Cauchy--Kovalevskaya product, used in
\cite{assv}, \cite{asv-cras} and \cite{MR2124899} to define
rational hyper-holomorphic functions, as we have discussed in the
previous section. The Wick product of two elements in the white
noise space  need not be in the white noise space; for an example
due to Gjessing, see \cite[Example 2.4.8 p. 45]{MR1408433}. There
is therefore a need to go beyond the white noise space.
Appropriate settings are Kondratiev's space $S_{-1}$ and Hida's
space of distributions. These spaces are defined below (see also
\cite[pp. 35-36]{MR1408433}), where the notation
\[
(2{\mathbb N})^{-q\alpha}=\prod_j (2j)^{-q\alpha_j}\]
is used.

\begin{Dn}
The Kondratiev space $S_{-1}$ consists of all formal power series
\eqref{fps} such that
\begin{equation}
\label{michelle} \sum_{\alpha}c_\alpha^2 (2{\mathbb
N})^{-q\alpha}<\infty
\end{equation}
for some $q\in{\mathbb N}$.\\
The Hida space $S^*$ consists of formal power series
\eqref{fps}
such that
\begin{equation}
\label{michelle1} \sup_{\alpha}c_\alpha^2 (2{\mathbb
N})^{-q\alpha}<\infty
\end{equation}
for some $q\in{\mathbb N}$.
\end{Dn}
We now introduce two spaces which  are the counterparts of the
Hardy space of the polydisk and of the Arveson space,
respectively. The case of the Arveson space is considered in
Section \ref{arv} and the polydisk case in Section \ref{280108}.

\begin{Dn}
The space ${\mathcal P}$ consists of all formal power series
\eqref{fps} for which
\begin{equation}
\sum_{\alpha\in\ell} c_\alpha^2<\infty.
\label{pld}
\end{equation}
The Arveson space ${\mathcal A}$ consists of all formal power
series \eqref{fps} for which
\begin{equation}
\sum_{\alpha\in\ell}\frac{\alpha!}{|\alpha|!}c_\alpha^2<\infty.
\label{arveson}
\end{equation}
\label{180108}
\end{Dn}

It is clear that ${\mathcal P}$ is included in the Hida space. We
will show in Theorem \ref{Jane_Birkin} below that ${\mathcal A}$
is included in the Kondratiev space.\\

The notion of the Wick product has been used in the development of
an Ito-like calculus for the fractional Brownian; see
\cite{bosw}, \cite{MR1741154}. The main properties of the Wick
product are listed in \cite[p. 43]{MR1408433}. We note in
particular the following:
\begin{enumerate}
\item The Wick product is independent of the basis; see
\cite[Appendix D, p. 209]{MR1408433}.

\item The Wick product differs in general from the point-wise
product. They coincide when at least one of the factors is
deterministic (see \cite[Example 2.4.6, p.43]{MR1408433}).

\item
The Wick product is not local. See \cite[2.4.10, p.45]{MR1408433}.
\end{enumerate}
A key property of the basis $(H_\alpha)$ is the following: define
a map ${\bf I}$ such that
\[
{\bf I}(H_\alpha)=z^\alpha,
\]
where $\alpha=(\alpha_1,\alpha_2,\ldots)\in\ell$, where
$z=(z_1,z_2,\ldots)\in {\mathbb C}^{\mathbb N}$ and where we use
the classical multi-index notation
\[
z^\alpha=z_1^{\alpha_1}\cdots
\]
Then,
\[
{\bf I}(H_\alpha\lozenge H_\beta)={\bf I}(H_\alpha){\bf
I}(H_\beta)=z^{\alpha+\beta}
\]
The map ${\bf I}$ is called the Hermite transform; it exhibits an
isomorphism between the white noise probability space and a
certain reproducing kernel Hilbert space, namely the space of
powers series
\[f(z)=\sum_{\alpha}z^\alpha f_\alpha,\quad\mbox{{\rm with norm}}\quad
\|f\|^2=\sum_{\alpha\in\ell}\frac{|f_\alpha|^2}{\alpha!}.\] This
space, called the Fock space, has been studied for a long time,
see for instance the 1962 paper of V. Bargmann \cite{bargmann}.
It is the reproducing kernel Hilbert space with reproducing kernel
\begin{equation}
K(z,w)=\sum_{\alpha\in\ell}
\frac{z^\alpha\overline{w}^\alpha}{\alpha!}= e^{\langle
z,w\rangle_{{\ell_2}}},\quad z,w\in{\ell_2}.
\end{equation}
See \cite[(10), p. 201]{bargmann}.\\

The analogue of the Cauchy-Kovalesvkaya extension theorem is the
following result; see \cite[Theorem 2.6.11, p. 62]{MR1408433}. In
the statement, $({\mathbb C}^{\mathbb N})_c$ denotes the space of
finite sequences of complex numbers indexed by the integers, and
(see \cite[Definition 2.6.4 p. 59]{MR1408433})
\[
K_q(\delta)=\left\{ z\in{\mathbb C}^{\mathbb
  N}\,\,:\,\,\sum_{\alpha\not=0}
|z|^\alpha
(2{\mathbb N})^{q\alpha}<\delta^2\right\}.
\]
\begin{Tm}\cite{MR1408433}
\label{ayefut} Let $g(z)=\sum_{\alpha\in\ell}g_\alpha z^\alpha$
be a power series defined in $({\mathbb C}^{\mathbb N})_c$, and
assume that $g$ is absolutely convergent in a domain $K_q^\delta$
for some $q<\infty$ and $\delta>0$. Then $g$ is the Hermite
transform of the element $G=\sum_{\alpha\in\ell} g_\alpha
H_\alpha$, which belongs to $S_{-1}$.
\end{Tm}

For $F=(F_{\ell j})\in (S_{-1})^{p\times q}$ we define ${\bf
I}(F)=({\bf I}(F_{\ell j}))$. The following propositions will be
used in the next section. The first claim is \cite[Proposition
2.2.6 p. 59]{MR1408433}, when $F$ is a row-valued function and
$G$ is a column-valued function. The result still holds for
matrices of appropriate dimensions, as is seen by taking component by
component. A similar remark holds for the second claim, which is
a consequence of Kondratiev's theorem (Theorem \ref{ayefut}
above), and corresponds to the function $f(x)=1-x$ in
\cite[Definition 2.6.14 p. 65]{MR1408433}.

\begin{Pn} \cite{MR1408433}
Let $F$ and $G$ be two matrix-valued functions with entries in
the Kondratiev space. Then
\[
{\bf I}(FG)={\bf I}(F){\bf I}(G)\quad{\rm and}\quad {\bf
I}(F+G)={\bf I}(F)+{\bf I}(G),
\]
where in each case, $F$ and $G$ are assumed to
be of compatible dimensions.
\end{Pn}
\begin{Pn}
\label{inv5}
For $F\in{S_{-1}}^{p\times p}$
\[
F=\sum_\alpha H_\alpha F_\alpha,\quad
F_\alpha\in{\mathbb R}^{p\times p}
\]
such that the constant coefficient is
$F_0=0_{p\times p}$, the von Neumann series
\[
(I-F)^{-\lozenge}=\sum_{k=0}^\infty F^{k\lozenge}
\]
converges in the Kontratiev space and
\[
{\bf I}((I_{p\times p}-F)^{-\lozenge})=(I_{p\times p}-{\bf
I}(F))^{-1}.
\]
\end{Pn}

Finally we recall that the Hermite transform allows to reduce
convergence in the Kondratiev space into convergence in terms of
power series. The following result holds (see \cite[Theorem 2.8.1
p. 74]{MR1408433}).

\begin{Tm} \cite{MR1408433}
A sequence of elements $F^{(n)}$ in the Kondratiev space
$S_{-1}$ converges to $F\in S_{-1}$ if there exists
$\delta>0$ and $q<\infty$ such that ${\bf I}(F^{(n)})$ converges
to ${\bf I}(F)$ pointwise boundedly, or equivalently, uniformly,
in $K_q(\delta)$.
\end{Tm}

As a corollary we give a direct application of this theorem.
Recall first that a Hilbert space ${\mathcal H}$ is called a
sub-Hilbert space of a topological vector space ${\mathcal V}$ if
it is included in ${\mathcal V}$ and if, moreover, the inclusion is
continuous. See \cite{schwartz}.

\begin{Cy}
\label{sofsof} Let ${\mathcal H}$ be a reproducing kernel Hilbert
space of functions $K(z,w)$ defined in a neighborhood
$K_q(\delta)$ and assume that $K(z,z)$ is uniformly bounded
there. Then, ${\mathcal H}$ is a subHilbert space of $S_{-1}$.
\end{Cy}

{\bf Proof:} It suffices to show that the inclusion map is continuous. Let
$(f_k)$ be a sequence of elements of ${\mathcal H}$ converging to
$f$ in the topology of ${\mathcal H}$. The reproducing kernel
property and the hypothesis on $K(z,z)$ implies that the pointwise
convergence is uniform in $K_q(\delta)$.\mbox{}\qed\mbox{}\\

\section{Rational functions}
\label{Anaelle}
\setcounter{equation}{0}
We now define and characterize rational functions in the white
noise space. We follow the methodology of the previous section. As
we have already mentioned, in the hyper-holomorphic case, the
Cauchy-Kovalevskaya theorem ensures that every ${\mathcal H}$--valued
function real analytic in a neighborhood of the origin, leads to a
hyper-holomorphic function. This result applies in particular when
the function is rational. In the stochastic setting, we will
define rational functions in terms of Hermite transforms. The
Cauchy-Kovalevskaya theorem is now replaced by Kondratiev's
theorem (Theorem \ref{ayefut} above), and we need to show
that a rational function satisfies the hypothesis of Theorem
\ref{ayefut}. This is done in the first theorem of this section.
In the statement we have set
\begin{equation}
\label{Hk} H_k=H_{e_k}
\end{equation}
where $e_k=(0,0,\ldots , 1,0,0,\ldots)$ is the element of $\ell$
with all entries $0$, at the exception of the $k$-th entry, which
is equal to $1$.

\begin{Tm}
Let $f$ be a ${\mathbb C}^{p\times q}$-valued rational function, analytic in a
neighborhood of the origin. Then $f$ is the image under the Hermite
transform of an element $F\in(S_{-1})^{p\times q}$. If
\[
f(z)=D+C(I_N-\sum_{k=1}^M z_kA_k)^{-1}(\sum_{k=1}^M z_kB_k)
\]
is a realization of $f$, then
\begin{equation}
\label{rational}
F=D+C(I_N-\sum_{k=1}^M
H_kA_k)^{-\lozenge}\lozenge (\sum_{k=1}^M H_kB_k).
\end{equation}
\end{Tm}
{\bf Proof:} The function $f$ is analytic in a neighborhood of the
closed ball
\[
B_{k,\epsilon}=\left\{(z_1,\ldots, z_M)\in{\mathbb C}^M\,\,;\,\,
\sum_{k=1}^M|z_k|^2\le \epsilon^2\right\}
\]
for some $\epsilon>0$. Fix $r_0>0$. We claim that for $q$ large
enough, it holds that
\begin{equation}
K_q(R_0)\subset B_{k,\epsilon}.
\label{incl}
\end{equation}
Indeed, let $z\in K_q(r_0)$. Then
\[
\sum_{\alpha\not=0}|z|^\alpha(2{\mathbb N}^{q\alpha})<r_0^2.\] In
particular,
\[
\sum_{k=1}^M |z_k|^{2}(2k)^{4qk}<r_0^2,\] and so
\[
\sum_{k=1}^M |z_k|^2<\frac{r_0^2}{2^{4q}}.
\]
For $q$ large enough, $\frac{r_0^2}{2^{4q}}\le \epsilon^2$. The
first claim of the theorem follows then from Theorem \ref{ayefut},
and the second claim follows from Proposition \ref{inv5}.
\mbox{}\qed\mbox{}\\

We now present the counterpart of Theorem \ref{tototototot}.
First a remark: The range of the Hermite functions consists of
functions which depend on a countable number of variables. We say
that such a function is {\it rational} if it depends only on a
finite number of these variables and is, moreover, a rational
function of these variables.

\begin{Tm}
Let $F\in (S_{-1})^{p\times q}$. The following are
equivalent:
\begin{enumerate}
\item ${\bf I}(F)$ is a rational function, analytic at the origin.
\item There are $N,M\in{\mathbb N}$ and matrices
$D\in{\mathbb R}^{p\times q}$, $C\in{\mathbb
R}^{p\times N}$, $A_1,\ldots ,A_M\in{\mathbb
R}^{N\times N}$ and $B_1,\ldots ,B_M\in{\mathbb
R}^{N\times q}$ such that $F$ is of the form \eqref{rational}.
\item $F$ is obtained from the $(H_\alpha)$  after
a finite number of additions, Wick multiplications and inversions
(the latter defined as in Proposition \ref{inv5}).
\end{enumerate}
\label{odeon}
\end{Tm}

{\bf Proof:} Assume that ${\bf I}(F)$ is rational and analytic in a
neighborhood of the origin. By the previous theorem, $F$ is then
of the asserted form. The converse follows by applying the Hermite
transform on $F$. The equivalence with the third statement follows
from the fact that a rational function is obtained after a finite
number of sums, products
and divisions of monomials\mbox{}\qed\mbox{}\\
\begin{Dn}
A function in the white noise space will be
called {\rm rational} if any of the equivalent
conditions in Theorem \ref{odeon} holds.
\end{Dn}

\begin{Tm} The Wick product and sum of two rational
functions of compatible dimensions are rational. If $F$ is rational and
invertible at the origin, then $F^{-1}$ is also rational.
\end{Tm}

{\bf Proof:} This follows from the first definition, since the
corresponding properties for rational functions
hold.\mbox{}\qed\mbox{}\\

\section{Realizable functions}
We now widen the class of rational functions.
\begin{Dn}
An element $F=\sum_\alpha c_\alpha H_\alpha$ in the white noise
space will be called {\rm realizable} if its Hermite transform can
be written as
\[
{\bf I}(F)(z)=D+C(I-zA)^{-1}zB
\]
where $D={\bf I}(F)(0)$, $z=(z_1,z_2,\ldots)$ and
\[
A=\begin{pmatrix}A_1\\ A_2\\
\vdots\end{pmatrix},\quad
B=\begin{pmatrix}B_1\\ B_2\\
\vdots\end{pmatrix},
\]
and where, moreover, the $A_j$ are operators
acting on a common Hilbert space ${ \mathcal H}$
and the $B_j$ are bounded operators from
${\mathbb C}$ into ${\mathcal H}$. We will say
that $F$ is {\rm finitely realizable} if ${\mathcal
H}$ is finite dimensional.
\end{Dn}

\begin{Ee}
Let $a\in\ell_2$ of norm strictly less than one, the function
\begin{equation}
\label{liora} b_a(z)=\frac{(1-|a|^2_{\ell_2})^{1/2}}{1-\langle z,a
\rangle_{\ell_2}}(z-a)(I_{\ell_2}-a^*a)^{-1/2}
\end{equation}
is realizable.
\end{Ee}

Indeed,
\[
b_a(z)=-(1-|a|^2_{\ell_2})^{1/2}a(I_{\ell_2}-a^*a)^{-1/2}+
(1-|a|^2_{\ell_2})^{1/2}(1-\langle z,a \rangle_{\ell_2})^{-1}
z(I_{\ell_2}-a^*a)^{1/2}.
\]
The function \eqref{liora} is the infinite dimensional version of
the Blaschke factors in the ball; it satisfies
\[
\frac{1-\langle b_a(z),b_a(w)\rangle_{\ell_2}}{1-\langle
z,w\rangle_{\ell_2}}= \frac{1-|a|^2_{\ell_2}}{(1-\langle
z,a\rangle_{\ell_2})(1-\langle a,w\rangle_{\ell_2})}.
\]
See \cite[p. 26]{rudin-ball}, \cite[pp. 11-13]{akap1} for a proof
in the unit ball of ${\mathbb C}^N$. The proof presented in this
last reference still holds in $\ell_2$.

\begin{Tm} The product and sum of two realizable
functions of compatible dimensions are realizable. If $F$ is realizable
and invertible at the origin, then $F^{-1}$ is also realizable.
\end{Tm}

{\bf Proof:} For non square matrices, addition is a special case
of multiplication since
\[
M+N=\begin{pmatrix}M&I_p\end{pmatrix}\begin{pmatrix}I_q\\
N\end{pmatrix},
\]
where $M$ and $N$ are $p\times q$ matrices. To check that the
product of two realizable functions is still realizable we
proceed as follows: let
\[
{\bf I}(F^{(k)})=D^{(k)}+C^{(k)}(I-zA^{(k)})^{-1}zB^{(k)},\quad
k=1,2
\]
be two realizable functions such that the product ${\bf
I}(F^{(1)}){\bf I}(F^{(2)})$ makes sense. Then,
\[
\begin{split}
{\bf I}(F^{(1)}){\bf
I}(F^{(2)})&=(D^{(1)}+C^{(1)}(I-zA^{(1)})^{-1}zB^{(1)})(
D^{(2)}+C^{(2)}(I-zA^{(2)})^{-1}zB^{(2)})\\
&=D^{(1)}D^{(2)}+\begin{pmatrix}C^{(1)}&D^{(1)}C^{(2)}\end{pmatrix}
\begin{pmatrix}I-zA^{(1)}&-zB^{(1)}C^{(2)}\\
0&I-zA^{(2)}\end{pmatrix}^{-1}\begin{pmatrix}zB^{(1)}D^{(2)}\\
zB^{(2)}\end{pmatrix}\\
&=D+C(I-zA)^{-1}zB,
\end{split}
\]
where $A,B,C$ and $D$ are defined as follows:
$A=\begin{pmatrix}A_1&A_2&\cdots\end{pmatrix}$ with
\[
A_\ell=\begin{pmatrix}A_\ell^{(1)}&B_\ell^{(1)}C^{(2)}\\0&A_\ell^{(2)}\end{pmatrix},\quad
\ell=1,2,\ldots,
\]
and
\[
B=\begin{pmatrix} B_1\\ B_2\\ \vdots\end{pmatrix},\quad{\rm
where}\quad
B_\ell=\begin{pmatrix} B_\ell^{(1)}D^{(2)}\\
B_\ell^{(2)}\end{pmatrix},\quad \ell=1,2,\ldots
\]
These formulas follow from
\[
\begin{split}
\begin{pmatrix}zA^{(1)}&zB^{(1)}C^{(2)}\\
0&zA^{(2)}\end{pmatrix}&=\begin{pmatrix}
z_1A_1^{(1)}+z_2A_2^{(1)}+\cdots&
z_1B_1^{(1)}C^{(2)}+z_2B_2^{(1)}C^{(2)}+\cdots\\
0&z_1A_1^{(2)}+z_2A_2^{(2)}+\cdots\end{pmatrix}
\\
&=z_1A_1+z_2A_2+\cdots,
\end{split}
\]
with the $A_\ell$ as above. A similar argument holds for $B$.
Moreover, $C$ and $D$ are given by the formula
\[
C=\begin{pmatrix}C^{(1)}&D^{(1)}C^{(2)}\end{pmatrix},\quad
D=D^{(1)}D^{(2)}.\]

The claim on the inverse follows from the formula
\[
(D+C(I-zA)^{-1}zB)^{-1}=D^{-1}-D^{-1}C(I-z(A-BD^{-1}C))^{-1}zBD^{-1}.
\]
This formula is easily verifiable through direct computation.
\mbox{}\qed\mbox{}\\
At present we do not have necessary and sufficient conditions
characterizing realizable functions. Two examples are provided in
the next two sections.

\section{The Arveson space and multipliers}
\label{arv}
We study here the Arveson space (see Definition \ref{180108}). Its
classical counterpart plays an important role in
multi-dimensional system theory. See \cite{btv}, \cite{MR1800795}.\\

\begin{Tm}
\label{Jane_Birkin}
The Arveson space is a subHilbert space of the
Kondratiev space.
\end{Tm}

{\bf Proof:} The Hermite transform of ${\mathcal A}$ is the
classical Arveson space, the reproducing kernel Hilbert space
with reproducing kernel
\begin{equation}
\label{k_ar}
k_{w}(z)=(1-\langle
z,w\rangle_{\ell_2})^{-1}=\sum_{\alpha\in\ell}
\frac{|\alpha|!}{\alpha!}z^\alpha\overline{w}^\alpha,
\end{equation}
where $z$ and $w$ are in the unit ball of $\ell_2$
\[
{\mathbb B}=\left\{z\in\ell_2\,\,;\,\,\sum_{k=1}^\infty
|z_k|^2<1\right\}.
\]
At this stage we wish to apply Corollary \ref{sofsof}. We note
that $K(z,z)$ is uniformly bounded in the closed ball of radius
$1/\sqrt{2}$
\[
\left\{z\in\ell_2\,\,;\,\,\sum_{k=1}^\infty |z_k|^2\le 1/2
\right\}.
\]
But this closed ball is in turn included in $K_q(1/\sqrt{2})$ for
any $q\ge 1$, and the proof is easily concluded.
\mbox{}\qed\mbox{}\\

\begin{Dn} A function $s$, defined in the open unit ball of $\ell_2$, is a
{\rm Schur multiplier} if, by definition, the operator of
multiplication by $s$ is a contraction from the classical Arveson
space into itself.
\end{Dn}

Equivalently, $s$ is a Schur multiplier if and only if the kernel
\[
K_s(z,w)=\frac{1-s(z)s(w)^*}{1-\langle z,w\rangle_{\ell_2}}
\]
is positive in the unit ball of $\ell_2$. We will denote by
${\mathcal H}(s)$ the corresponding reproducing kernel Hilbert
space and \( \Gamma_s=I-M_sM_s^*\), where $M_s$ is the operator
of multiplication in the Arveson space. We have
\[
K_s(z,w)=(\Gamma_s(k_w))(w), \]
and the range of $\Gamma_s$ is dense in ${\mathcal H}(s)$, as
follows from general results on operator ranges and reproducing
kernel spaces; see for instance \cite{fw}.
\begin{Tm}
A Schur multiplier is the image under the Hermite transform of a
realizable function in the Kondratiev space.
\end{Tm}

{\bf Proof:} A number of related proofs hold for this result in
the case where only a finite number of variables are involved;
see for instance \cite{btv}, \cite{adr1}. Here we briefly outline the
method of \cite[Th\'eor\`eme 2.1]{asv-cras}, suitably adapted to the case
of a countable number of variables. We consider the case of a
scalar-valued Schur multiplier (the general case is treated in the
same way). We proceed with the following steps:\\

STEP 1: {\sl The operators of multiplication by the variables
$M_{z_k}$ are bounded in ${\mathcal A}$ and it holds that
\begin{equation}
\label{fuitedeau}
\sum_{k=0}^\infty M_{z_k}M_{z_k}^*=I-C^*C,
\end{equation}
where $Cf=f(0)$.}\\

One first verifies that $M_{z_k}^*={\mathcal R}_k$. The identity
\eqref{fuitedeau} is then obvious.\\

STEP 2: {\sl Let ${\mathcal H}(s)_\infty$ be the closure in
${\mathcal H}(s)\oplus{\mathcal H}(s)\oplus \cdots$ of the
functions
\[
w_y=\begin{pmatrix}\Gamma_sM_{z_0}^*k_y\\ \Gamma_sM_{z_1}^*k_y\\
\vdots\end{pmatrix}
\]
where $y$ runs through the open unit ball of $\ell_2$, and where
$k_y$ is defined by \eqref{k_ar}. The formulas
\[
\widetilde{G}(1)=K_s(z,0),\quad \widetilde{H}(1)=s(0)^*
\]
and
\[
\widetilde{T}(w_y)=K_s(z,y)-K_s(z,0),\quad
\widetilde{F}(w_y)=s(y)^*-s(0)^*
\]
define an isometric relation on ${\mathcal H}(s)_\infty\oplus
{\mathbb C}$ into ${\mathcal
H}(s)\oplus{\mathbb C}$ with a dense domain.}\\

The proof is as in \cite{asv-cras}, and uses equality \eqref{fuitedeau}.\\

Since the linear relation in the previous step is isometric and
densely defined, it is the graph of an isometric operator. We
denote by
\[\begin{pmatrix} T&F\\ G&H\end{pmatrix}\]
its adjoint.\\


STEP 3: {\sl It holds that
\[
s(z)=H+G(I-\sum_{k=0}^\infty z_kT_k)^{-1}(\sum_{k=0}^\infty
z_kF_k).\] }

As in \cite{asv-cras}, one first proves that
\[
\begin{split}
\sum_{k=0}^\infty z_k(T_kf)(z)&=f(z)-f(0)\\
\sum_{k=0}^\infty z_k(F_kf)(z)&=s(z)-s(0).
\end{split}
\]
The formula for $s$ is then a direct consequence of these equations.
\mbox{}\qed\mbox{}\\

\section{The space ${\mathcal P}$ and Schur-Agler classes}
\setcounter{equation}{0}
\label{280108}
We now present another class of realizable functions, related to
the space ${\mathcal P}$. Recall that this space was defined
Definition \ref{180108}. As in the previous section, we consider
the scalar case to simplify notation.
\begin{Dn}
A function $s$ is in the Schur-Agler class if there exist a
family $k_1(z,w),k_2(z,w),\ldots$ of functions positive in
$\prod_{k=1}^\infty {\mathbb D}$, where ${\mathbb D}$ denotes the
open unit disk, and such that
\begin{equation}
\label{arc-de-triomphe}
1-s(z)s(w)^*=\sum_{\ell=1}^\infty (1-z_\ell w_\ell^*)k_\ell(z,w)
\end{equation}
there.
\end{Dn}
In the case of the finite polydisk, these classes originate with
the work of Agler and have been much studied; see
\cite{agler-hellinger}, \cite{MR2154356}, \cite{MR2336046}.

\begin{Tm}
A Schur-Agler multiplier is the Hermite transform of a realizable
element in the Kondratiev space.
\end{Tm}
{\bf Proof:} We follows the proof given in \cite[Section
4.2]{MR2240272}. Let ${\mathcal H}(k_\ell)$ be the reproducing
kernel Hilbert space of functions with reproducing kernel $k_\ell$
and let
\[
{\mathcal H}=\oplus_{\ell=1}^\infty {\mathcal H}(k_\ell),\]
with norm $(\sum_{\ell=1}^\infty\|f_\ell\|_{{\mathcal
H}(k_\ell)}^2)^{1/2}$. One defines a linear relation on
${\mathcal H}\oplus{\mathbb C}\times {\mathcal H}\oplus{\mathbb
C}$ via the formula
\[
\widehat{A}\begin{pmatrix} w_1^* k(\cdot, w_1)\\ w_2^* k(\cdot,
w_2)\\ \vdots\end{pmatrix}=\begin{pmatrix}
 w_1^* (k(\cdot, w_1)-k(\cdot, 0))\\ w_2^* (k(\cdot,
w_2)-k_2(\cdot, 0))\\ \vdots\end{pmatrix},\quad
\widehat{B}\begin{pmatrix} w_1^* k(\cdot, w_1)\\ w_2^* k(\cdot, w_2)\\
\vdots\end{pmatrix}=s(w)^*-s(0)^*,
\]
and
\[
\widehat{C}(1)=\begin{pmatrix} k_1(\cdot,0)\\ k_2(\cdot,0)\\
\vdots \end{pmatrix},\quad \widehat{D}=s(0)^*.
\]
Then
\[
\begin{pmatrix}\widehat{A}&\widehat{C}\\
\widehat{B}&\widehat{D}\end{pmatrix}
\]
extends to an isometric relation, which is the graph of an
operator. We denote by
\[
\begin{pmatrix}{A}&{B}\\
{C}&{D}\end{pmatrix}
\]
its adjoint. We have the formula
\[
z(Af)(z)=\sum_{\ell=0}^\infty (f_\ell(z)-f_\ell(0)),\quad
z(B(z))=s(z)-s(0),
\]
and
\[
Ch=\sum_{\ell=0}^\infty f_\ell(0),\quad D=s(0).
\]
We note that the infinite sums above converge because of the
Cauchy-Schwartz inequality. For instance
\[
\begin{split}
|\sum_{\ell=0}^\infty f_\ell(0)|&=|\sum_{\ell=0}^\infty\langle
f_\ell,
k_\ell(\cdot,0)\rangle_{{\mathcal H}(k_\ell)}\\
&\le \sum_{\ell=0}^\infty \|f_\ell\|_{{\mathcal
H}(k_\ell)}\sqrt{k_\ell(0,0)}\\
&\le(\sum_{\ell=0}^\infty\|f_\ell\|_{{\mathcal
H}(k_\ell)}^2)^{1/2}(\sum_{\ell=0}^\infty k_\ell(0,0))^{1/2}.
\end{split}
\]
Finally, one has the realization
$$s(z)=D+C(I-\sum_{\ell=1}^\infty
z_\ell A_\ell)^{-1}(\sum_{\ell=1}^\infty z_\ell B_\ell),
$$
where $A_\ell=\pi_\ell A$ and $B_\ell=\pi_\ell B$, with
$\pi_\ell$ being the orthogonal projection from ${\mathcal H}$
onto ${\mathcal H}_\ell$.
\mbox{}\qed\mbox{}\\
\section{The Gleason problem in certain spaces of power series}
\setcounter{equation}{0}
\label{soucot2007}
Consider a power series in a countable number of variables
$f(z)=\sum_{\alpha\in\ell}z^\alpha f_\alpha$, convergent in a
neighborhood of the origin. The Leibenzon's operators
${\mathcal R}_j$ can still be defined as power series as in
\eqref{bastille}, and we have:
\begin{equation}
\label{gleason} f(z)-f(0)=\sum_{j=0}^\infty z_j {\mathcal
R}_jf(z).
\end{equation}


We will say that a space of power series in the $z_j$ is
{\it resolvent-invariant} if Gleason's problem is solvable with bounded
operators $A_j$:
\begin{equation*}
f(z)-f(0)=\sum_{j=0}^\infty z_j (A_jf)(z).
\end{equation*}
The space will be called {\it backward-shift invariant} if the $A_j$
commute.\\

We now prove a uniqueness theorem in a collection of spaces of
functions depending on a countable number of variables, which
includes in particular the Fock space, the Arveson space and the
infinite polydisk space. The argument follows the one given in
the hyper-holomorphic setting in \cite{MR2240272}. For an earlier
result in the setting of power series in finite number of
variables, see \cite{adubi-pams}. Note that in the theorem we do
not assume the operator of multiplication by ${z_k}$ to be
bounded. The theorem holds in particular in the Arveson space and
in the Fock space. In the Fock space the operators of
multiplication by the variables are not bounded.

\begin{Tm}
\label{shirly} Let ${\mathcal H}$ be a reproducing kernel Hilbert
space of functions depending on a countable number of variables
$z_1,\ldots, $ defined in a neighborhood of the origin, and
assume that the Leibenzon operators ${\mathcal R}_j$ are
uniformly bounded in ${\mathcal H}$. Any other uniformly
commuting solution to Gleason's problem, coincides with the
${\mathcal R}_j$.
\end{Tm}

{\bf Proof:} Let $A_1,\ldots$ be a family of commuting and
uniformly bounded operators in ${\mathcal H}$ such that
\[
f(z)-f(0)=\sum_{n=0}^\infty z_n(A_nf)(z).
\]
Since the $A_j$ commute and are uniformly bounded, we can
reiterate this equation and obtain
\[
f(z)=\sum_{\alpha\in\ell} z^\alpha CA^\alpha f=C(I-zA)^{-1}f,\]
where $C$ is the operator of evaluation at the origin. In
particular,
\begin{equation}
\label{uniq} C(I-zA)^{-1}f\equiv 0\Longrightarrow f=0,
\end{equation}
and replacing $f$ by $A_kf$ we have:
\[
(A_kf)(z)=C(I-zA)^{-1}A_kf.
\]
On the other hand by definition of ${\mathcal R}_k$,
\[
({\mathcal R}_k
f)(z)=\sum_{\alpha\ge e_k}z^{\alpha-e_k}\frac{\alpha_k}{|\alpha|}CA^\alpha
f=C(I-zA)^{-1}A_kf,\]
and so $A_k={\mathcal R}_k$ in view of \eqref{uniq}.
\mbox{}\qed\mbox{}\\

As a corollary of the preceding theorem we have:
\begin{Tm}
The Leibenzon operators are bounded in the Fock space, in the
Arveson space and in the infinite polydisk space. In particular,
they are the only commutative solution of Gleason's problem in
these spaces.
\end{Tm}
{\bf Proof:} Let $f(z)=\sum_{\alpha\in\ell}z^\alpha f_\alpha$ be
in the Fock space ${\mathcal F}$. By definition of ${\mathcal
R}_j$ and of the norm in the Fock space, we have:
\[
\begin{split}
\|{\mathcal R}_jf\|_{\mathcal F}&=\sum_{\alpha\ge
e_j}\frac{\alpha_j^2}{|\alpha|^2}f_{\alpha}^2(\alpha-e_j)!\\
&\le \sum_{\alpha\ge e_j} f_\alpha^2\alpha!\\
&\le \|f\|^2_{\mathcal F}.
\end{split}
\]
\mbox{}\qed\mbox{}\\
%
\newpage
We now give the table presenting the parallels between the
hyper-holomorphic case and the stochastic case.\\

\begin{tabular}{|l|l|l|}
\hline & &\\
The setting & Hyper-holomorphic case & Stochastic case\\
& &\\
 \hline& &\\
The underlying space& Functions hyper-holomorphic & The
Kondratiev's space
$S_{-1}$\\
&at the origin&
\\
 \hline & &\\
 The building blocks& The
hyperholomorphic & The
functions $H_{e_k}$, $k\in{\mathbb N}$\\
&variables $\zeta_j$, $j=1,2,3$&\\
 \hline
Power series expansions
 &&\\
in terms of &
 the Fueter polynomials& the functions  $H_\alpha$\\
 \hline
 & &\\
 The product& Cauchy-Kovalesvkaya& Wick product\\
& &\\
 \hline
Commutativity& When the restrictions to& Always commutes\\
of the product& to $x_0=0$ commute&\\
 \hline Convolution on "power series"&  &\\
expansions &$f\circ g=\sum_{\alpha}\zeta^\alpha(\sum_{\beta\le
 \alpha}f_\beta g_{\alpha-\beta})$&
$f\lozenge g=\sum_{\alpha\in\ell}H_\alpha(\sum_{\beta\le
 \alpha}f_\beta g_{\alpha-\beta})$\\
 & &\\
 \hline
Going to the classical case from
 & Restriction to & The Hermite transform\\
the hyper-holomorphic/stochastic &the hyperplane $x_0=0$&\\
case& &\\
 \hline
From the classical case to the & &\\
hyper-holomorphic/stochastic
 &CK extension& Kontradiev's theorem\\
case& &(see Theorem \ref{ayefut})\\
 \hline
Uniqueness theorem for&In particular &In particular\\
the Leibenzon's operators& in the Arveson space& in
the Fock space\\
 \hline
\end{tabular}

\newpage
\bibliographystyle{plain}
\def\cprime{$'$} \def\lfhook#1{\setbox0=\hbox{#1}{\ooalign{\hidewidth
  \lower1.5ex\hbox{'}\hidewidth\crcr\unhbox0}}} \def\cprime{$'$}
  \def\cprime{$'$} \def\cprime{$'$} \def\cprime{$'$} \def\cprime{$'$}


\begin{thebibliography}{10}

\bibitem{agler-hellinger}
J.~Agler.
\newblock {\em On the representation of certain holomorphic functions defined
  on a polydisk}, volume~48 of {\em {Operator {T}heory: {A}dvances and
  {A}pplications}}, pages 47--66.
\newblock Birkh{\" a}user Verlag, Basel, 1990.

\bibitem{MR2275397}
D.~Alpay, F.~M. Correa-Romero, M.~E. Luna-Elizarrar{\'a}s, and M.~Shapiro.
\newblock Hyperholomorphic rational functions: the {C}lifford analysis case.
\newblock {\em Complex Var. Elliptic Equ.}, 52(1):59--78, 2007.

\bibitem{adr1}
D.~Alpay, A.~Dijksma, and J.~Rovnyak.
\newblock A theorem of {B}eurling--{L}ax type for {H}ilbert spaces of functions
  analytic in the ball.
\newblock {\em Integral Equations Operator Theory}, 47:251--274, 2003.

\bibitem{ad4}
D.~Alpay and C.~Dubi.
\newblock {Backward shift operator and finite dimensional de Branges Rovnyak
  spaces in the ball}.
\newblock {\em {Linear Algebra and Applications}}, 371:277--285, 2003.

\bibitem{a-dubi3}
D.~Alpay and C.~Dubi.
\newblock A realization theorem for rational functions of several complex
  variables.
\newblock {\em {System \& Control Letters}}, 49:225--229, 2003.

\bibitem{adubi-pams}
D.~Alpay and C.~Dubi.
\newblock On commuting operators solving {G}leason's problem.
\newblock {\em Proc. Amer. Math. Soc.}, 133(11):3285--3293 (electronic), 2005.

\bibitem{akap1}
D.~Alpay and H.T. Kaptano\u{g}lu.
\newblock Some finite-dimensional backward shift-invariant subspaces in the
  ball and a related interpolation problem.
\newblock {\em Integral Equation and Operator Theory}, 42:1--21, 2002.

\bibitem{assv}
D.~Alpay, B.~Schneider, M.~Shapiro, and D.~Volok.
\newblock Fonctions rationnelles et th\'eorie de la r\'ealisation: le cas
  hyper--analytique.
\newblock {\em {C}omptes {R}endus {M}ath\'ematiques}, 336:975--980, 2003.

\bibitem{as1}
D.~Alpay and M.~Shapiro.
\newblock Gleason's problem and tangential homogeneous interpolation for
  hyperholomorphic quaternionic functions.
\newblock {\em {Complex Variables}}, 48:877--894, 2003.

\bibitem{asv-cras}
D.~Alpay, M.~Shapiro, and D.~Volok.
\newblock Espaces de de {B}ranges {R}ovnyak: le cas hyper--analytique.
\newblock {\em {C}omptes {R}endus {M}ath\'ematiques}, 338:437--442, 2004.

\bibitem{MR2124899}
D.~Alpay, M.~Shapiro, and D.~Volok.
\newblock Rational hyperholomorphic functions in {$R^4$}.
\newblock {\em J. Funct. Anal.}, 221(1):122--149, 2005.

\bibitem{MR2240272}
D.~Alpay, M.~Shapiro, and D.~Volok.
\newblock Reproducing kernel spaces of series of {F}ueter polynomials.
\newblock In {\em Operator theory in Krein spaces and nonlinear eigenvalue
  problems}, volume 162 of {\em Oper. Theory Adv. Appl.}, pages 19--45.
  Birkh\"auser, Basel, 2006.

\bibitem{MR2154356}
J.~Ball and V.~Bolotnikov.
\newblock Nevanlinna-{P}ick interpolation for {S}chur-{A}gler class functions
  on domains with matrix polynomial defining function in {$\Bbb C\sp n$}.
\newblock {\em New York J. Math.}, 11:247--290 (electronic), 2005.

\bibitem{MR2336046}
J.~Ball and V.~Bolotnikov.
\newblock Interpolation in the noncommutative {S}chur-{A}gler class.
\newblock {\em J. Operator Theory}, 58(1):83--126, 2007.

\bibitem{btv}
J.~Ball, T.~Trent, and V.~Vinnikov.
\newblock Interpolation and commutant lifting for multipliers on reproducing
  kernel {H}ilbert spaces.
\newblock In {\em Proceedings of Conference in honor of the 60--th birthday of
  {M.A. Kaashoek}}, volume 122 of {\em Operator {T}heory: {A}dvances and
  {A}pplications}, pages 89--138. Birkhauser, 2001.

\bibitem{bargmann}
V.~Bargmann.
\newblock Remarks on a {H}ilbert space of analytic functions.
\newblock {\em Proceedings of the {N}ational {A}cademy of {Arts}}, 48:199--204,
  1962.

\bibitem{bgk1}
H.~Bart, I.~Gohberg, and M.A. Kaashoek.
\newblock {\em Minimal factorization of matrix and operator functions},
  volume~1 of {\em {Operator {T}heory: {A}dvances and {A}pplications}}.
\newblock Birkh{\" a}user Verlag, Basel, 1979.

\bibitem{bosw}
F.~Biagini, B.~{\O}ksendal, A.~Sulem, and N.~Wallner.
\newblock An introduction to white-noise theory and {M}alliavin calculus for
  fractional {B}rownian motion, stochastic analysis with applications to
  mathematical finance.
\newblock {\em Proc. R. Soc. Lond. Ser. A Math. Phys. Eng. Sci.},
  460(2041):347--372, 2004.

\bibitem{MR924203}
M.~Bisiacco, E.~Fornasini, and G.~Marchesini.
\newblock Controller design for {$2$}{D} systems.
\newblock In {\em Frequency domain and state space methods for linear systems
  (Stockholm, 1985)}, pages 99--113. North-Holland, Amsterdam, 1986.

\bibitem{bds}
F.~Brackx, R.~Delanghe, and F.~Sommen.
\newblock {\em Clifford analysis}, volume~76.
\newblock Pitman research notes, 1982.

\bibitem{MR1741154}
T.~Duncan, Y.~Hu, and B.~Pasik-Duncan.
\newblock Stochastic calculus for fractional {B}rownian motion. {I}. {T}heory.
\newblock {\em SIAM J. Control Optim.}, 38(2):582--612 (electronic), 2000.

\bibitem{eh}
R.J. Elliott and J.~van~der Hoek.
\newblock A general fractional white noise theory and applications to finance.
\newblock {\em Math. Finance}, 13(2):301--330, 2003.

\bibitem{fw}
P.A. Fillmore and J.P. Williams.
\newblock On operator ranges.
\newblock {\em {Advances in Mathematics}}, 7:254--281, 1971.

\bibitem{MR54:2342}
E.~Fornasini and G.~Marchesini.
\newblock State--space realization theory of two--dimensional filters.
\newblock {\em IEEE Trans. Automatic Control}, AC--21(4):484--492, 1976.

\bibitem{MR80c:93028}
E.~Fornasini and G.~Marchesini.
\newblock Doubly-indexed dynamical systems: state-space models and structural
  properties.
\newblock {\em Math. Systems Theory}, 12(1):59--72, 1978/79.

\bibitem{MR35:7123}
I.M. Guelfand and N.Y. Vilenkin.
\newblock {\em Les distributions. {T}ome 4: {A}pplications de l'analyse
  harmonique}.
\newblock Collection Universitaire de Math\'ematiques, No. 23. Dunod, Paris,
  1967.

\bibitem{MR1408433}
H.~Holden, B.~{\O}ksendal, J.~Ub{\o}e, and T.~Zhang.
\newblock {\em Stochastic partial differential equations}.
\newblock Probability and its Applications. Birkh\"auser Boston Inc., Boston,
  MA, 1996.

\bibitem{MR93j:26013}
S.~G. Krantz and H.~P. Parks.
\newblock {\em A primer of real analytic functions}, volume~4 of {\em Basler
  Lehrb\"ucher [Basel Textbooks]}.
\newblock Birkh\"auser Verlag, Basel, 1992.

\bibitem{MR1800795}
S.~McCullough and T.~Trent.
\newblock Invariant subspaces and {N}evanlinna-{P}ick kernels.
\newblock {\em J. Funct. Anal.}, 178(1):226--249, 2000.

\bibitem{rudin-ball}
W.~Rudin.
\newblock {\em Function theory in the unit ball of ${\mathbb{C}}^n$}.
\newblock Springer--{V}erlag, 1980.

\bibitem{Rudin}
W.~Rudin.
\newblock {\em Real and complex analysis}.
\newblock {Mc Graw Hill}, 1982.

\bibitem{schwartz}
L.~Schwartz.
\newblock Sous espaces hilbertiens d'espaces vectoriels topologiques et noyaux
  associ\'{e}s (noyaux reproduisants).
\newblock {\em J. Analyse Math.}, 13:115--256, 1964.

\bibitem{MR618518}
F.~Sommen.
\newblock A product and an exponential function in hypercomplex function
  theory.
\newblock {\em Applicable Anal.}, 12(1):13--26, 1981.

\bibitem{MR54:11066}
E.L. Stout.
\newblock {\em The theory of uniform algebras}.
\newblock Bogden \& Quigley, Inc., Tarrytown-on-Hudson, N. Y., 1971.

\end{thebibliography}
Daniel Alpay\\
Department of Mathematics\\
Ben--Gurion University of the Negev\\
Beer-Sheva 84105, Israel\\
{\tt dany@math.bgu.ac.il}\\

David Levanony\\
Department of Electrical Engineering\\
Ben--Gurion University of the Negev\\
Beer-Sheva 84105, Israel\\
{\tt levanony@ee.bgu.ac.il}\\
\end{document}